\documentclass[11pt]{article}
\usepackage{euscript,amsmath,amssymb }
\usepackage{srcltx}

\title{On two classes of dense 2-generator  subgroups in $\mathbb C$.}

\author{Kirill Kamalutdinov \and 
Andrey Tetenov\footnote{Supported by Russian Foundation of Basic Research project
13-01-00513} \and Dmitry Vaulin}

\begin{document}

\newcommand{\rr}{\mathbb{R}}
\newcommand{\qq}{\mathbb{Q}}
\newcommand \nn {\mathbb{N}}
\newcommand \zz {\mathbb{Z}}
\newcommand \bbc {\mathbb{C}}
\newcommand \rd {\mathbb{R}^d}

 \newcommand {\al} {\alpha}
\newcommand {\be} {\beta}
\newcommand {\da} {\delta}
\newcommand {\Da} {\Delta}
\newcommand {\ga} {\gamma}
\newcommand {\Ga} {\Gamma}
\newcommand {\la} {\lambda}
\newcommand {\La} {\Lambda}
\newcommand{\om}{\omega}
\newcommand{\Om}{\Omega}
\newcommand {\sa} {\sigma}
\newcommand {\Sa} {\Sigma}
\newcommand {\te} {\theta}
\newcommand {\fy} {\varphi}
\newcommand {\ep} {\varepsilon}
\newcommand{\e}{\varepsilon}

\newcommand{\VEC}{\overrightarrow}
\newcommand{\IN}{{\subset}}
\newcommand{\NI}{{\supset}}
\newcommand \dd  {\partial}
\newcommand {\mmm}{{\setminus}}
\newcommand{\probel}{\vspace{.5cm}}
\newcommand{\8}{{\infty}}
\newcommand{\0}{{\varnothing}}
\newcommand{\vse}{$\blacksquare$}

\newcommand {\bfep} {{{\bar \varepsilon}}}
\newcommand {\Dl} {\Delta}
\newcommand{\vA}{{\vec {A}}}
\newcommand{\vB}{{\vec {B}}}
\newcommand{\vF}{{\vec {F}}}
\newcommand{\vf}{{\vec {f}}}
\newcommand{\vh}{{\vec {h}}}
\newcommand{\vJ}{{\vec {J}}}
\newcommand{\vK}{{\vec {K}}}
\newcommand{\vP}{{\vec {P}}}
\newcommand{\vX}{{\vec {X}}}
\newcommand{\vY}{{\vec {Y}}}
\newcommand{\vZ}{{\vec {Z}}}
\newcommand{\vx}{{\vec {x}}}
\newcommand{\va}{{\vec {a}}}
\newcommand{\vga}{{\vec {\gamma}}}

\newcommand{\eS}{{\EuScript S}}
\newcommand{\eH}{{\EuScript H}}
\newcommand{\eC}{{\EuScript C}}
\newcommand{\eP}{{\EuScript P}}
\newcommand{\eT}{{\EuScript T}}
\newcommand{\eG}{{\EuScript G}}
\newcommand{\eK}{{\EuScript K}}
\newcommand{\eF}{{\EuScript F}}
\newcommand{\eZ}{{\EuScript Z}}
\newcommand{\eL}{{\EuScript L}}
\newcommand{\eD}{{\EuScript D}}
\newcommand{\E}{{\EuScript E}}
\def \diam {\mathop{\rm diam}\nolimits}
\def \fix {\mathop{\rm fix}\nolimits}
\def \Lip {\mathop{\rm Lip}\nolimits}

\newtheorem{thm}{\bf Theorem}
 \newtheorem{cor}[thm]{\bf Corollary}
 \newtheorem{lem}[thm]{\bf Lemma}
 \newtheorem{prop}[thm]{\bf Proposition}
 \newtheorem{dfn}[thm]{\bf Definition}

\newcommand{\dok}{{\bf{Proof}}}

\maketitle

\begin{abstract} We consider dense 2-generator multiplicative  subgroups in $\mathbb C$ and show that for each point $z\in\bbc $ the set of limit values for the arguments of the powers of  each generator at the point $z$ is either finite or is $[-\pi,\pi]$
\end{abstract}
\maketitle
\vspace{12mm}\\

\section*{Introduction.}

The paper was motivated by some neat examples in fractal geometry \cite{ATK, KTV}  and the problems arising in the course of their construction.
\bigskip

It is a well known and  widely used  fact  that if  $\xi,\eta$ are such positive numbers,  that $\dfrac{\log \xi}{\log \eta}$ is irrational, then for  any $\la>0$ we can find such sequence
$(n_k,m_k)$, that  $\lim\limits_{k\to\8}n_k=\8$,$\lim\limits_{k\to\8} m_k=\8$, $\lim\limits_{k\to\8} \xi^{n_k}/\eta^{m_k}=\la$. 

\bigskip

Now suppose that $\xi$ and $\eta$ are complex numbers. What conditions must  be imposed on $\xi$  and $\eta$ to ensure  that such sequence  exists for any $\la\in\bbc$?  More complicated  question is, can we find such $\xi,\eta$, that for any $\la\in\bbc$ and for any $\al\in  (-\pi,\pi)$ there is a sequence $(n_k,m_k)$,
for which  $\lim\limits_{k\to\8} \xi^{n_k}/\eta^{m_k}=\la$, while the values of the arguments of $\xi^{n_k}, \eta^{m_k}$ converge to a given value $\al$?

Surprisingly, the answer is yes. This  is possible in the case, when $a,b$ are the generators of  such group  $G=\langle\xi,\eta,\cdot\rangle$,  dense in $\bbc$, that for  each open set $V\IN \bbc$ the values of  the  arguments of first factors $\xi^m$ of the  elements $\xi^m\eta^n\in (G\cap V)$ are  dense in $[-\pi,\pi]$.  We  call such $G$  {\em a group of  the second type} and prove  that  they do exist.

\section*{ Dense additive  and multiplicative groups in $\bbc$.}

We begin with dense 3-generator  lattices in $\bbc$.

\begin{thm} Let  $\xi=\al +i \be$, $\al,\be \in \rr$. A group $ G=\langle 1,i,\xi,+ \rangle $ is  dense in $\bbc$ iff
 for any integers $k,l,m$ 
$$ k\al+l\be+m =0  \mbox{\rm \ \  implies \ \ }  k=l=m=0, \eqno{(*)}$$
\end{thm}

{\bf Proof.}  

Let $\fy$ be the canonical homomorphism of $ G$  to the factor-group $G'= G/{\zz[i]} $, where $\zz[i]=\langle 1,i,+ \rangle$ is  a group of  Gaussian integers. $G'$ is a subgroup of the torus $T=\bbc/{\zz[i]}   $. By (*), the equality $\fy( m\xi+k+li)=\fy( m'\xi+k'+l'i)$ holds iff $m=m'$. Therefore all  $z_m=\fy( m\xi+k+li) $ are different in $G'$. So $G'$ is an infinite cyclic subgroup in $T$ and it has a limit point in $T$. If a sequence $\{z_{m_k}\}$ converges, the  sequence 
$\{z_{m_{k+1}}-z_{m_k}\}$ converges to 0, so 0 is   a limit  point  both in $G'$ and in   $ G$.

\bigskip

For  $z\neq 0$ in $\bbc$, let $L_z$ be the line $\{tz \   |  t\in\rr\}.$ 
Notice, that  for  any bounded  subset $D\IN \bbc$, the set $ G\cap L_z\cap D$ is finite.

 Indeed, if $z=m\xi+k+li,$ and $ z'=m'\xi+k'+l'i\in L_z$, then
$\dfrac{m'\al+k'}{m\al+k}=\dfrac{m'\be+l'}{m\be+l}$.  The last equation implies
$$(m'l-l'm)\al+(km'-k'm)\be+ (k'l-kl')=0  $$
Then it follows  from (*) that 
$\dfrac{l'}{l}=\dfrac{m'}{m}=\dfrac{k'}{k}$. If $z'\in D$, this  ratio is bounded, so $ G\cap D\cap L_z$ is finite.    

\bigskip

Therefore, for any $\e>0$  we can take such  $a,b \in  G\cap B(0,\e)$, that $b\notin L_a$. Then the  set $ \{ k a+l b :\ \  k,l\in \zz \}$ is  an $\e$-net in $\bbc$.  This shows that $ G$ is  dense in $\bbc$. 

\bigskip
The  second  part  of  the proof is rather short.

Suppose  $k\al+l\be+m =0 $ and $k\neq 0$.  Find such $p_1,p_2,q\in \zz$ that
$kp_1+lp_2-mq=0$. If $x=q\al+p_1$,  $y=q\be+p_2$ then $kx+ly=kp_1+lp_2-qm$. 

The last equality means that all the elements of the group $ G$ lie on a family of parallel lines $kx+ly=n$, where $k,l$ are fixed and $n\in \zz$. The unit  square P intersects no more  than $k+l+1$ of these lines, therefore $ G$ is not  dense in $\bbc$.\vse

Applying  an affine  tranformation sending $1,i,\xi $ to $u,v,w$, we  get

\begin{cor}\ Let $u,v \in \bbc$,  $Im \dfrac{u}{v} \neq 0$, $w=\al u + \be v$, $\al,\be \in \rr$. \medskip \\
A group $ G=\langle u,v,w,+ \rangle $ is  dense in $\bbc$ iff
 for any integers $k,l,m$, 
$$ k\al+l\be+m =0  \mbox{\rm \ \  implies \ \ }  k=l=m=0 \ \ \blacksquare$$

\end{cor}

In the  case $w=1$ we can apply  the map $f(z)=e^{2\pi i z}$ to obtain 

\begin{cor}\label{mulgr} Let $u,v \in \bbc$,  $Im \dfrac{u}{v} \neq 0$, $\al u + \be v=1$, $\al,\be \in \rr$, $\xi=e^{2\pi i u}$,  \medskip \\ $\eta=e^{2\pi i v}$.
A group $G=\langle \xi,\eta,\cdot\rangle $ is  dense in $\bbc$
 iff
 for any integers $k,l,m$, 
$$ k\al+l\be+m =0  \mbox{\rm \ \  implies \ \ }  k=l=m=0   \ \  \blacksquare$$
\end{cor} 

Then $G=\langle \xi,\eta,\cdot\rangle $ is called a  {\em dense 2-generator multiplicative  group}  in $\bbc$.

We will consider a group $G=\langle \xi,\eta,\cdot\rangle $ along with its additive counterpart, $\hat G=\langle u,v,1,+\rangle$, where $\xi=e^{2\pi i u},  \eta=e^{2\pi i v}$.

\bigskip

Notice that if  $G=\langle \xi,\eta,\cdot\rangle $ is dense in $\bbc$, the  formula $\psi(\xi^m\eta^n)=\left(\dfrac{\xi}{|\xi|}\right)^m$ \smallskip \\  defines 
  a homomorphism $\psi$ of  a group $G=\langle \xi,\eta,\cdot\rangle $ to the unit circle
$S^1\IN \bbc$.
Put $$H_G= \bigcap\limits_{\e>0}\overline{\psi(B(1,\e)\cap G)}$$

In other  words, $H_G$ is the  set of limit points of all those  sequences $\{e^{i n_k\arg(\xi)}\}$, for  which   $\{n_k\}$ is the first  coordinate  projection  of such   sequence $\{(n_k,m_k)\}$, that $ \lim\limits_{k\to\8} \xi^{n_k}\eta^{m_k}=1$.

The set $H_G$ is a closed  subset of the unit  circle $S^1$, and it is  a subgroup of the  group $S^1$, so it is either finite cyclic, i.e. $H_G=\{e^{2k\pi i /n}\}$ for  some $n\in \nn$, or it is infinite and  dense in $S^1$ and  therefore $H_G=S^1$.

\begin{dfn}
A dense 2-generator subgroup $G$ is called  the group of the  first  type, if $H_G$ is finite, and  the group of the  second type, if $H_G=S^1$.
\end{dfn}

\begin{lem}\label{alph}  $G$ is of the  second type iff for  some $\al\notin \qq$,  $e^{2\pi i \al}\in H_G$.\vse
\end{lem}

{\bf Remark.} Let $\xi=e^{2\pi i u},  \eta=e^{2\pi i v}$ are  the  generators  of  the group $G$ of  the  second  type and $e^{2\pi i \al}\in H_G$ for some irrational $\al$. From the point of view of the group $G'$ this means, that there is such sequence $\{(p_n,r_n,s_n)\}$ in $\zz^3$, that $\lim\limits_{n\to\8}(r_{n}u+s_n v) =1$  and $\lim\limits_{n\to\8} r_{n}Re(u)-p_n =\al$. 

 In other words,
$\lim\limits_{n\to\8}\left\{ r_{n}Re(u)\right\} =\al$ (where $\{x\} $ stands for  fractional part of x).

\begin{cor}
If   $G=\langle \xi,\eta,\cdot\rangle$ is  of the  second type, then for  any  open subset $V\in \bbc$ the  set 
$\psi(V\cap G) $ is  dense in $S^1$.
\end{cor}

\dok \ \ Take $\xi^m\eta^n\in V$. Then for some $\e>0$, $ \xi^{-m}\eta^{-n}V\NI B(1,\e)$. Since $\psi(B(1,\e)$ is dense in $S^1$, the same is true for $\psi(V)$. \vse

The  groups of  the  second type  have  a significant  geometric  property:

\begin{thm}
If  the  group $G$ with  generators $\xi=r e^{i\al},\eta=R e^{i\be}$ is  of the  second type, then for  any   $z_1,z_2\in \bbc\mmm\{0\}$ there is such sequence $\{(n_k,m_k)\}$  that $\lim\limits_{k\to \8}   \dfrac{z_1\xi^{n_k}}{z_2\eta^{m_k}}=1$  ,  $\lim\limits_{k\to\8} e^{i n_k\al}=e^{-i \arg(z_1)}$ ,  $\lim\limits_{k\to\8} e^{i n_k\be}=e^{-i \arg(z_2)}$.
\end{thm}

The proof  directly follows from the  previous corollary.\vse

\bigskip

The  last  theorem means,  that if $G$ is  of the  second type,  then given any  two orbits  $\{z_1\xi^{n}\},\{z_2\eta^{m}\}$ and   any ray $L_\theta:\arg z=\theta$, there are two subsequences $\{z_1\xi^{n_k}\},\{z_2\eta^{m_k}\} $ for  which $\lim\limits_{k\to \8}   \dfrac{z_1\xi^{n_k}}{z_2\eta^{m_k}}=1$, while both projections
$\left(\dfrac{z_1\xi^{n_k}}{|z_1\xi^{n_k}|}\right),\left(\dfrac{z_2\eta^{m_k}}{|z_2\eta^{m_k}|}\right) $ approach $e^{i\theta}$.

% For the group $\hat G=\langle u,v,1,+\rangle$, corresponding to the  group $G$ of  the  second type, 

\section*{  Examples of  the groups of first  and second  type.}

For the  groups of  the  first  type, we have  the following  sufficient  condition:
\begin{thm}
Suppose $G=\langle \xi,\eta,\cdot\rangle$ is  dense  in $\bbc$, $\xi=e^{2\pi iu}$,
$\eta=e^{2\pi iv}$. 
If  $$\dfrac{{\rm Re}u\,{\rm Im}v}{{\rm Re}v\,{\rm Im}u}\in\mathbb{Q},$$
then $G$ is of the  first type.

\end{thm}
{\bf Proof.}
Put $u=u_{1}+iu_{2}$, $v=v_{1}+iv_{2}$. Suppose $\dfrac {v_{1}u_{2}}{u_{1}v_{2}}=\dfrac{p}{q}$,
where $p,\, q\in\zz$. Notice that $u$ and $v$ are not collinear, so $p\neq q$.

Take some $\e\in(0,1)$ and suppose $|\xi^m\eta^n-1|<\e$, $m,n\in\zz$.
Then $|mu+nv-k|<\e$  for some $k\in\zz$.\\
Combining the inequalities $|nu_{2}+mv_{2}|<\e$,
 $|mu_{1}+nv_{1}-k|<\e$,
we get \bigskip \\ $\left|m\dfrac{u_{1}v_{2}-v_{1}u_{2}}{v_{2}}-k\right|<\e\left(1+\left|\dfrac{v_{1}}{v_{2}}\right|\right)$, or
 $\left| m\dfrac{q-p}{q}u_{1}-k\right|<\e\left(1+\left|\dfrac{v_{1}}{v_{2}}\right|\right). $ \bigskip

Rewriting  the  last  inequality in the  form $$\left|mu_1-\dfrac{kq}{q-p}\right|<\e\left|\dfrac{q}{q-p}\right|\left(1+\left|\dfrac{v_{1}}{v_{2}}\right|\right)$$ we  see  that  the order of the  group $G$ divides $|p-q|$.\vse

\bigskip

To show the  existence of  a dense  group of  the  second  type, we follow  Remark to the  Lemma 5 and find  such $u,v\in \bbc$ and $\al\in (0,1)\mmm \qq$, that
for  some integers $r_n, s_n$,  $\lim\limits_{n\to\8}(r_{n}u+s_n v) =1$ and $\lim\limits_{n\to\8}\left\{ r_{n}Re(u)\right\} =\al$. 
We make it  three  steps.

{\bf  Step 1.}  Construct irrational numbers $\al, \be, \ga$.

  We construct such positive numbers  $\alpha$,   $\be$ and  $\ga$ with continued  fraction representations
$$\al=\cfrac{1}{a_{1}+\cfrac{1}{a_{2}+...}},\quad \quad \be=\cfrac{1}{b_{1}+\cfrac{1}{b_{2}+...}},\quad \quad  \ga=\cfrac{1}{c_{1}+\cfrac{1}{c_{2}+...}},$$ that  the denominators of the
convergents for  $\al, \be, \ga$ form the sequences $\{q_n\}, \ \  \{q_n+1\}, \ \  \{q_n +2 + (-1)^{n}\}$ respectively.
So the  convergents for  $\al, \be, \ga$ will be
$$\dfrac{p_n}{q_n}, \ \  \dfrac{r_n}{q_n+1}, \ \  \dfrac{s_n}{q_n +2 + (-1)^{n}}$$

\bigskip

Define the sequence $\{q_{n}\}=(1,2,31,994,...)$  by the relations: $q_{1}=1$, $q_{2}=2$,
$q_{2n+1}=q_{2n-1}+q_{2n}(q_{2n}+1)(q_{2n}+3)$,
$q_{2n+2}=q_{2n}+q_{2n+1}(q_{2n+1}+1)$.

\bigskip

The terms  for $\alpha$
are (see \cite{K}) \ \ \  $a_{n+2}=\dfrac{q_{n+2}-q_{n}}{q_{n+1}}$,\medskip \\
so $a_{2n+1}=(q_{2n}+1)(q_{2n}+3)$,
$a_{2n+2}=q_{2n+1}+1$, $a_{1}=q_{1}=1$. 

\bigskip

The terms for $\be$ are $b_{2n+1}=q_{2n}(q_{2n}+3)$,
$b_{2n+2}=q_{2n+1}$, $b_{1}=2$.

\bigskip

 Then the terms for $\ga$ are  $c_{2n+1}=q_{2n}(q_{2n}+1)$,
$c_{2n+2}=q_{2n+1}$,  $c_{1}=2$.

\bigskip

Finally we get the following sequences $a_n, b_n, c_n$:\\ $(1,1,15,32, 995\cdot 997,...),$\\
$(2,1,10,31, 994\cdot 997,...),$\\  $ (2,2,6,32, 994\cdot 995,...)$.\\

\bigskip

For odd n, we have
 $\be=\dfrac{r_{n}}{q_{n}+1}+\dfrac{\be_{n}}{q_{n}+1}$ and
$\ga=\dfrac{s_{n}}{q_{n}+1}+\dfrac{\ga_{n}}{q_{n}+1}$. \medskip \\
Then $\be_{n}<\dfrac{1}{q_{n}+1}$ and 
$\ga_{n}<\dfrac{1}{q_{n}+1}$ for any $n\in \nn$.

\bigskip

{\bf  Step 2.}  Show that $\be$ and $\ga$ satisfy the condition (*) of Theorem 1.\bigskip

Suppose $l\be+m\ga-k=0$ for some  $l,m,k\in\mathbb{Z}$.

Take such $N$ that for $n>N$,  $q_{n}>3\max\left\{  |k |, |m |, |l |\right\} $.

Rewrite  $l\be+m\ga-k=lr_{n}+ms_{n}-k(q_{n}+1)+l\be_{n}+m\ga_{n}=0$.
Since $l\be_{n}+m\ga_{n}<\dfrac{l+m}{q_{n}+1}<\dfrac{2}{3}$, we have $\be_{n}=\ga_{n}=0$ for $n>N$. Then for odd $n>N$, we  get $$l\dfrac{r_{n}}{q_{n}+1}+m\dfrac{s_{n}}{q_{n}+1}=k  \mbox{\ \ and \ \ \ }
 l\dfrac{r_{n+1}}{q_{n+1}+1}+m\dfrac{s_{n+1}}{q_{n+1}+3}=k $$

Therefore $l\left(\dfrac{r_{n}}{q_{n}+1}-\dfrac{r_{n+1}}{q_{n+1}+1}\right)+m\left(\dfrac{s_{n}}{q_{n}+1}-\dfrac{s_{n+1}}{q_{n+1}+3}\right)=0$.

\bigskip

Observing that $r_{n}(q_{n+1}+1)-r_{n+1}(q_{n}+1)=1$
and $s_{n}(q_{n+1}+3)-s_{n+1}(q_{n}+1)=1$,
we get  that for any odd $n>N$, $(l+m)q_{n+1}+3l+m=0$. This is possible only if $l=m=0$, so $k$ is also 0.

\bigskip

{\bf Step 3.}    Find  such $u,v\in\bbc$  that
$\be u+\ga v=1$, $Re(\be u)=\alpha$ to construct the group $G$ . 

To get desired $u,v$, take $h>0$ and  define $u=\dfrac{\al +ih}{\be}$, $v=\dfrac{1-\al -ih}{\ga}$.

For any $\e>0$ there is   such  $N$, that for  any   $n>N$,   $|q_{n}\al-p_{n}|<\e$,
$|(q_{n}+1)\be-r_{n}|<\e$  and $|(q_{n}+2+(-1)^n)\ga-s_{n}|<\e$. \medskip

For any odd $n>N$,  the sum of  the second and third inequalities, multiplied by $|u|,|v|$ respectively, 
gives  $|(q_{n}+1)-(r_{n}u +s_{n}v)|<\e(|u| +|v| )$.

From $ |(q_{n}+1)\al-r_{n}Re (u) |<\e|u|$
and $ |(q_{n}+1)\al-(p_{n}+\al) |<\e$
we get $ |r_{n}Re(u)-p_{n}-\alpha |<\e(|u| +1)$. \bigskip

Using that $0<\al<1$, and $p_{n}$ is an integer,  for $n=2k+1$ we have 
$\lim\limits_{k\to\8}\left\{ r_{2k+1}Re(u)\right\} =\al$ (where $\{x\} $ stands for  fractional part of x), while
$\lim\limits_{k\to\8}(r_{2k+1}u+s_{2k+1} v) =1$ .

Consider  the  group $G=\langle \xi,\eta,\cdot\rangle$ generated by $\xi=e^{2i\pi u}$ and $\eta=e^{2i\pi v}$.
By  Corollary \ref{mulgr}, G is dense in $\bbc$. Since 
$$\lim\limits_{k\to\8}e^{ i r_{2k+1}Re(u)} =e^{2\pi i\al}, $$  where $\al$ is irrational,
by Lemma \ref{alph}  the group  G is of the second type.


\begin{thebibliography}{99}

\bibitem{ATK} Aseev V.~V.,Tetenov A.~V.,Kravchenko A.~S., On Self-Similar Jordan Curves on the Plane, Siberian Mathematical Journal,  May/Jun  2003, Vol.~44, Issue~3, p.~379


\bibitem{K} Khinchin A.~Ya., Continued fractions / Mineola, N.Y.: Dover Publications, 1997.

\bibitem{KTV} Kamalutdinov K.G., Tetenov A.V., Vaulin D.A., Self-Similar  Jordan Arcs, which do not  satisfy Weak Separation Property, (in prepapation)


\end{thebibliography}
\end{document}